\documentclass{amsart}
\usepackage[T1]{fontenc}

\usepackage{graphicx,amsmath,amssymb,amsfonts,bm,float,cite}
\usepackage[shortlabels]{enumitem}
\usepackage{kpfonts}
\usepackage[none]{hyphenat}
\usepackage[hyphens]{url}
\usepackage[hypertexnames=false]{hyperref}
\usepackage[center]{caption}
\usepackage{mathtools}

\DeclarePairedDelimiter{\ceil}{\lceil}{\rceil}
\DeclarePairedDelimiter{\floor}{\lfloor}{\rfloor}

\makeatletter
\edef\orig@output{\the\output}
\output{\setbox\@cclv\vbox{\unvbox\@cclv\vspace{0pt plus 20pt}}\orig@output}

\title{The rectangular spiral \\
or the $n_1 \times n_2 \times \cdots \times n_k$ Points Problem}
\author{Marco Rip\`a}

\begin{document}

\begin{abstract}
\sloppy A generalization of Rip\`a's square spiral solution for the $n \times n \times \cdots \times n$ Points Upper Bound Problem. Additionally, we provide a non-trivial lower bound for the $k$-dimensional $n_1 \times n_2 \times \cdots \times n_k$ Points Problem. In this way, we can build a range in which, with certainty, all the best possible solutions to the problem we are considering will fall. Finally, we give a few characteristic numerical examples in order to appreciate the fineness of the result arising from the particular approach we have chosen.
\end{abstract}

\maketitle


\section{Introduction}

\sloppy Nearly a century ago, the classic \emph{Nine Dots Problem} appeared in Sam Loyd's \emph{Cyclopedia of $5000$ Puzzles, Tricks and Conundrums with Answers} \cite{1, 4}. The challenge was as follows: ``\ldots draw a continuous line through the center of all the eggs so as to mark them off in the fewest number of strokes'' \cite{3, 5}.

\begin{figure}[H]
\includegraphics[width=\columnwidth, keepaspectratio]{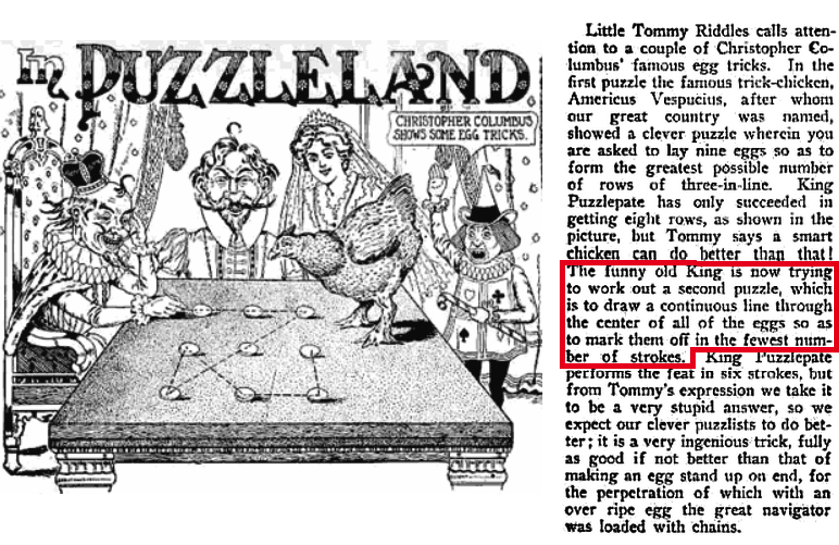}
\caption{The original problem published in Sam Loyd's \emph{Cyclopedia of $5000$ puzzles}, New York, 1914, p. 301.}
\label{fig:Figure_1R}
\end{figure}

That puzzle can naturally be extended to an arbitrarily large number of distinct (zero-dimensional) points for each row/column \cite{6}. This new problem asks to connect $n \times n$ points, arranged in a grid formed by $n$ rows and $n$ columns, using the fewest straight lines connected at their endpoints. Rip\`a and Remirez \cite{6} showed that it is possible to do this for every $n \in \mathbb{N} - \{0,1,2\}$, using only $2 \cdot n-2$ straight lines. For each $n \geq 5$, we can combine a given $8$-line solution for the $5 \times 5$ problem and the square spiral frame \cite{10}. In the same paper, they extended the $n \times n$ result to a three-dimensional space \cite{7} and finally to a generic $k$-dimensional space (for $k>3$).

Starting from that outcome, we consider the same problem and rules as Reference \cite{6}. We can apply the ``pure'' spiral method to a $n_1 \times n_2$ rectangular grid (where $n_1 \leq n_2$). In this way, it is simple to discover that the minimum number of line segments connected at their endpoints required to join every point (solving the problem inside the box, connecting dots without crossing a line, and visiting each point just once) cannot exceed $h_u$, the (planar) square spiral upper bound given by (\ref{eq1}) \cite{9},
\begin{equation} \label{eq1}
h_u = 2 \cdot n-1,  \forall n \in \mathbb{N}-\{0,1\}.
\end{equation}


\section{The \texorpdfstring{$n_1 \times n_2 \times \cdots \times n_k$}{} Problem upper bound}

If we try to extend the result in (\ref{eq1}) to a three-dimensional space, where \linebreak$n_1 \leq n_2 \leq n_3$, we need to modify somewhat the standard strategy described in Reference \cite{6} in order to choose a good ``plane-by-plane'' approach that we can implement. So, we need to identify the correct starting plane to lay the first straight line. 
Using basic mathematics, it is quite easy to prove that, in general, the best option is to start from the $[n_2; n_3]$ plane.

Hence, under the additional constraints that we must solve the problem inside the box only, connecting points without crossing a line, and visiting each dot just once, our strategy is as follows:
\begin{itemize}
    \item 	Step 1) Take one of the external planes identified by $[n_2; n_3]$, and here is the plane to lay our first line segment;
    \item 	Step 2) Starting from one corner of this planar grid, draw the first straight line to connect $n_3$ points until we reach the last point in that row;
    \item 	Step 3) The second line segment is on the same plane. It lays on $[n_2; n_3]$, it is orthogonal to the previous one, and it links $n_2-1$ points;
    \item 	Step 4) Repeat the square (rectangular) spiral pattern until we connect all the $n_2 \cdot n_3$ points of the mentioned planar grid to the others on the same surface;
    \item 	\sloppy Step 5) Draw the ($n_2-1$)-th line, which is orthogonal to the $[n_2; n_3]$ plane we have considered before, visiting $n_1-1$ new points by crossing all the remaining planar grids, and then double the same scheme (in reverse) on the opposite face of this parallelogram (inside the three-dimensional axis-aligned bounding box $[n_1 \times n_2 \times n_3]$). Repeat this pattern for each $n_2 \times n_3$ grid, $n_1-2$ more times.
\end{itemize}

The rectangular spiral solution also gives us the shortest path that visits every point of any given $n_1 \times n_2 \times n_3$ grid since the total Euclidean length of the line segments used to fit all the points is minimal.
\vspace{4mm}

\sloppy \noindent {\bf Nota Bene.} Just a couple of trivial considerations. Referring to the rectangular spiral pattern applied to a $k$-dimensional space ($k \geq 2$), we can return to the starting point using exactly one additional line (it works for any number of dimensions). For any odd value of $n_1$, we can visit a maximum of $\ceil[\big]{\frac{n_{k-1}}{2}}-1$ points twice, simply extending the line end (if we do not, we will not visit any dot more than once, otherwise we can visit $\ceil[\big]{\frac{n_{k-1}}{2}}-1$ points, at most). Moreover, it is possible to visit up to $n_{k-1}-2$ points twice if we move the second to last line too (crossing some more lines as well). Finally, assuming $k \geq 2$, if we are free to extend the ending line until we are close to the next (already visited) point (i.e., let $\varepsilon$ be the distance between the last line and the nearest point and let the distance between two adjacent points be unitary, we have that $0< \varepsilon <1$), it is possible to return to the starting point without visiting any point more than once.

The number of line segments we spend to connect every point is always lower than or equal to
\begin{equation} \label{eq2}
h_u \leq 2 \cdot n_1 \cdot n_2-1.
\end{equation}

In fact, $h_u \leq (2 \cdot n_2-1) \cdot n_1  + n_1-1$.

Nevertheless, $(2 \cdot n_2-1)\cdot n_1+n_1-1=2 \cdot n_1 \cdot n_2-n_1+n_1-1=2 \cdot n_1 \cdot n_2-1=2 \cdot n_1 \cdot n_2-n_2+n_2-1=(2 \cdot n_1-1) \cdot n_2+n_2-1$.  (Q.E.D.)

The ``savings'', in terms of unused segments, are zero if (and only if)
\begin{equation} \label{eq3}
n_1 < 2 \cdot (n_3-n_2)+3.
\end{equation}

In general, (also if $n_1 \geq 2 \cdot (n_3-n_2)+3$), (\ref{eq2}) can be improved as
\begin{equation} \label{eq4}
h_u=2 \cdot n_1 \cdot n_2-c,
\end{equation}
where $c=1$ if the savings are zero, while $c \geq 2$ if not.

As an example, let us consider the following cases:
\begin{enumerate}
\item $n_1=5$, $n_2=6$, $n_3=9$;
\item $n_1=11$, $n_2=12$, $n_3=13$.
\end{enumerate}
While in the first instance $c=1$ (since $5 < 2 \cdot (9-3)+3)$ so that $h=2 \cdot 5 \cdot 6-1=59$, in case b) we have $c=13$, and thus $h=2 \cdot 11 \cdot 12-13=251$. This is by virtue of the fact that the fifth and the sixth connecting lines allow us to ``save'' one line for every subsequent plane, whereas each plane ``met'' after the sixth can be solved using two fewer lines (if compared with the first four we have considered).

\begin{figure}[H]
\includegraphics[scale=0.24, keepaspectratio]{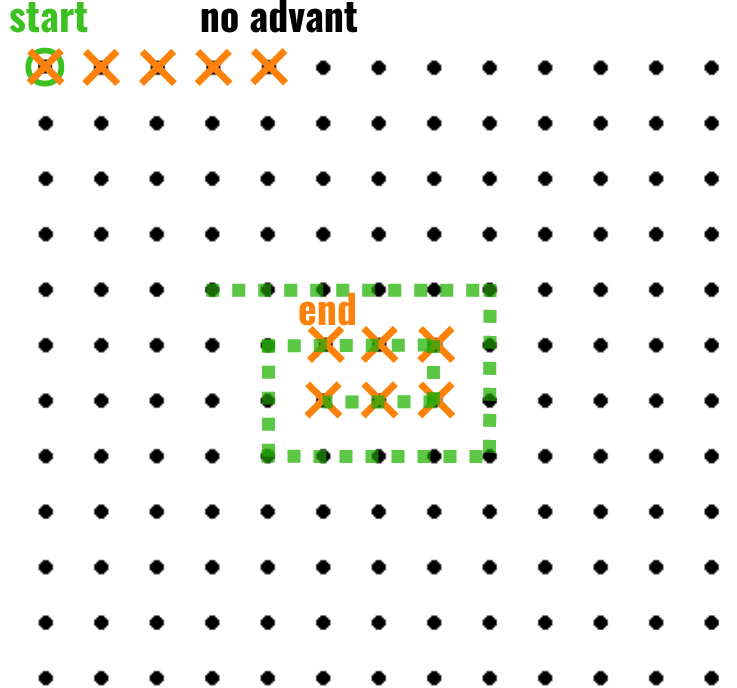}
\caption{The rectangular spiral for the case of the example b): $n_1=11$, $n_2=12$, $n_3=13$.}
\label{fig:Figure_2S}
\end{figure}

If $n_1 \geq 2 \cdot (n_3-n_2)+4$, the (pure) rectangular spiral method, with specific regard to the three-dimensional problem, can be summarized as
\begin{gather*}
h_u=n_1-1+[2 \cdot (n_3-n_2)+2] \cdot (2 \cdot n_2-1)+2 \cdot (2 \cdot n_2-2)+[2 \cdot (n_3-n_2)+4] \cdot \\ (2 \cdot n_2-3)+ 4 \cdot (2 \cdot n_2-4)+[2 \cdot (n_3-n_2)+6] \cdot (2 \cdot n_2-5)+6 \cdot (2 \cdot n_2-6)+ \cdots + d,
\end{gather*}
where $d$ represents the product of the number of line segments used to solve the plane which contains the fewest lines (the last planar grid we have considered, laying on the plane which cuts roughly halfway through our imaginary \emph{box}) and the residual $n_1-\{[2 \cdot (n_3-n_2)+2]+2+[2 \cdot (n_3-n_2)+4]+4+\cdots \}$.

Then, let $b$ and $j_{max}$ be non-negative integers to synthesize the formula above as
$$h_u=n_1-1+\sum_{j=0}^{j_{max}} [(2 \cdot n_2-2 \cdot j-1) \cdot(2 \cdot(n_3-n_2)+2 \cdot (j+1))+2 \cdot(j+1) \cdot(2 \cdot n_2-2 \cdot(j+1))]+b.
$$

Hence,
\begin{equation} \nonumber
\begin{aligned}
h_u= -\frac{8 \! \cdot \! {j_{max}}^3}{3}+6 \! \cdot \! {j_{max}}^2 \! \cdot \! n_2 - 2 \! \cdot \! {j_{max}}^2  \! \cdot n_3 - 11 \! \cdot \! {j_{max}}^2 - 4 \! \cdot \! j_{max} \! \cdot \! {n_2}^2 + 4  \! \cdot \! j_{max}  \! \cdot n_2 \! \cdot \! n_3 + \\ 16 \! \cdot \! j_{max}  \! \cdot \! n_2 - 4 \! \cdot j_{max} \! \cdot \! n_3 - \frac{43 \! \cdot \! j_{max}}{3} - 4 \! \cdot \! {n_2}^2 + 4 \! \cdot \! n_2 \! \cdot \! n_3 + 10 \! \cdot \! n_2 - 2 \! \cdot \! n_3 - 7 + n_1 + b,
\end{aligned}
\end{equation}
where $j_{max}$ represents the \underline{maximum} value of the upper bound of the summation, let us say $\tilde{j}$, such that 
$$
n_1 \geq \sum_{j=0}^{\tilde{j}} [2 \cdot (n_3-n_2)+2 \cdot (j+1)+2 \cdot (j+1)] \Rightarrow n_1 \geq 2 \cdot (\tilde{j}+1) \cdot (n_3-n_2+\tilde{j}+2),
$$
while
\begin{equation} \nonumber
\begin{aligned}
b \coloneqq 
    \begin{cases}
[n_1- 2 \cdot (j_{max}+1) \cdot (n_3-n_2+j_{max}+2)] \cdot (2 \cdot n_2-2 \cdot j_{max}-3)  \\ \mathbf{if} \hspace{5mm} n_1-2 \cdot (j_{max}+1) \cdot (n_3-n_2+j_{max}+2) \leq 2 \cdot (n_3-n_2)+2 \cdot (j_{max}+2)  \vspace{5mm} \\
[2 \cdot (n_3-n_2)+2 \cdot (j_{max}+2)] \cdot (2 \cdot n_2-2 \cdot j_{max}-3)+ \{n_1-2 \cdot (j_{max}+1) \cdot \\ (n_3-n_2+j_{max}+2)- [2 \cdot (n_3-n_2)+2 \cdot (j_{max}+2)]\} \cdot (2 \cdot n_2-2 \cdot j_{max}-4) \\ \mathbf{if} \hspace{5mm} n_1- 2 \cdot (j_{max}+1) \cdot (n_3-n_2+j_{max}+2) > 2 \cdot (n_3-n_2)+2 \cdot (j_{max}+2)
\end{cases}.
\end{aligned}
\end{equation}

By making some calculations, we have that
\begin{equation} \nonumber
\begin{aligned}
b \coloneqq 
    \begin{cases}
4 \cdot {j_{max}}^3-8 \cdot {j_{max}}^2 \cdot n_2+4 \cdot {j_{max}}^2 \cdot n_3+18 \cdot {j_{max}}^2-2 \cdot j_{max} \cdot n_1+4 \cdot j_{max} \cdot {n_2}^2- \\ 4 \cdot j_{max}  \cdot n_2 \cdot n_3-22 \cdot j_{max}  \cdot n_2+10 \cdot j_{max}  \cdot n_3+26 \cdot j_{max}+2 \cdot n_1 \cdot n_2-3 \cdot n_1+ \\ 4 \cdot {n_2}^2-4 \cdot n_2 \cdot n_3-14 \cdot n_2+6 \cdot n_3+12
\\ \mathbf{if} \hspace{5mm} n_1 \leq 2 \cdot (j_{max}+2) \cdot (j_{max}-n_2+n_3+2) \vspace{5mm}  \\
4 \cdot {j_{max}}^3-8 \cdot {j_{max}}^2 \cdot n_2+4 \cdot {j_{max}}^2 \cdot n_3+20 \cdot {j_{max}}^2-2 \cdot j_{max} \cdot n_1+4 \cdot j_{max} \cdot {n_2}^2- \\ 4 \cdot j_{max}  \cdot n_2 \cdot n_3-24 \cdot j_{max}  \cdot n_2+12 \cdot j_{max}  \cdot n_3+34 \cdot j_{max}+2 \cdot n_1 \cdot n_2-4 \cdot n_1+ \\ 4 \cdot {n_2}^2-4 \cdot n_2 \cdot n_3-18 \cdot n_2+10 \cdot n_3+20 \\ \mathbf{if} \hspace{5mm} n_1>2 \cdot (j_{max}+2) \cdot (j_{max}-n_2+n_3+2)
    \end{cases}\hspace{-1mm}.
\end{aligned}
\end{equation}

Thus, the general solution is given by:
\begin{equation} \label{eq5}
\begin{aligned}
h_u =
    \begin{cases}
\frac{4 \cdot {j_{max}}^3}{3}-2 \cdot {j_{max}}^2 \cdot n_2+2 \cdot {j_{max}}^2 \cdot n_3+7 \cdot {j_{max}}^2-2 \cdot j_{max} \cdot  n_1-6 \cdot j_{max} \cdot  n_2 + \\ 6 \cdot j_{max} \cdot  n_3+ \frac{35 \cdot j_{max}}{3}+2 \cdot n_1 \cdot n_2-2 \cdot n_1-4 \cdot n_2+4 \cdot n_3+5
\\ \mathbf{if} \hspace{5mm} n_1 \leq 2 \cdot ({j_{max}}^2-j_{max} \cdot  n_2+j_{max} \cdot n_3+4 \cdot j_{max}-2 \cdot n_2+2 \cdot n_3+4)  \vspace{5mm}  \\
\frac{4 \cdot {j_{max}}^3}{3}-2 \cdot {j_{max}}^2 \cdot n_2+2 \cdot {j_{max}}^2 \cdot n_3+9 \cdot {j_{max}}^2-2 \cdot j_{max} \cdot  n_1-8 \cdot j_{max} \cdot  n_2+ \\ 8 \cdot j_{max} \cdot  n_3+ \frac{59 \cdot j_{max}}{3}+2 \cdot n_1 \cdot n_2-3 \cdot n_1-8 \cdot n_2+8 \cdot n_3+13 \\ \mathbf{if} \hspace{5mm} n_1 > 2 \cdot ({j_{max}}^2-j_{max} \cdot  n_2+j_{max} \cdot n_3+4 \cdot j_{max}-2 \cdot n_2+2 \cdot n_3+4)
    \end{cases},
\end{aligned}
\end{equation}
where $j_{max}$ is the \underline{maximum} non-negative integer $j$ such that $$n_1 \geq 2 \cdot [j^2+(n_3-n_2+3) \cdot j+n_3-n_2+2],$$ so $j_{max}= \floor[\big]{ \frac{1}{2} \cdot \left(\sqrt{{n_3}^2+{n_2}^2-2 \cdot n_2 \cdot n_3+2 \cdot n_3-2 \cdot n_2+2 \cdot n_1+1}+n_2-n_3-3 \right)}$.

Then, (\ref{eq5}) can be more elegantly written as
\begin{equation} \label{eq6}
\begin{aligned}
h_u =
    \begin{cases}
\frac{4}{3}  \cdot {j_{max}}^3+[2 \cdot (n_3-n_2)+7]  \cdot {j_{max}}^2+\left[6  \cdot (n_3-n_2)-2  \cdot n_1+\frac{35}{3} \right]  \cdot j_{max}+ \\ 4 \cdot (n_3-n_2)+2 \cdot n_1 \cdot (n_2-1)+5
\\ \mathbf{if} \hspace{5mm} n_1 \leq 2 \cdot [{j_{max}}^2+ (n_3-n_2+4) \cdot j_{max}+2 \cdot (n_3-n_2)+4]  \vspace{5mm}  \\
\frac{4}{3}  \cdot {j_{max}}^3+[2 \cdot (n_3-n_2)+9]  \cdot {j_{max}}^2+\left[8 \cdot (n_3-n_2)-2 \cdot  n_1+\frac{59}{3} \right]  \cdot j_{max}+ \\ 8 \cdot (n_3-n_2 )+n_1 \cdot (2 \cdot n_2-3)+13 \\ \mathbf{if} \hspace{5mm} n_1>2 \cdot [{j_{max}}^2+ (n_3-n_2+4) \cdot j_{max}+2 \cdot (n_3-n_2)+4]
    \end{cases},
\end{aligned}
\end{equation}
where $j_{max}=\floor[\Big]{\frac{\sqrt{(n_3-n_2+1)^2+2 \cdot n_1}-n_3+n_2-3}{2}}$.
\vspace{4mm}

\noindent {\bf Nota Bene.} For obvious reasons, (\ref{eq6}) is always applicable, on condition that \linebreak$n_1 \geq 2 \cdot (n_3-n_2)+4$. Otherwise, the solution follows immediately from (\ref{eq4}) since $c$ can assume only two distinct values: $1$ or $2$ ($c=1$ if (\ref{eq3}) is verified, whereas $c=2$ if (\ref{eq3}) is not satisfied and (\ref{eq6}) cannot be used -- therefore, this is the case $n_1 = 2 \cdot (n_3-n_2)+3$).

\begin{figure}[H]
\includegraphics[width=\columnwidth, keepaspectratio]{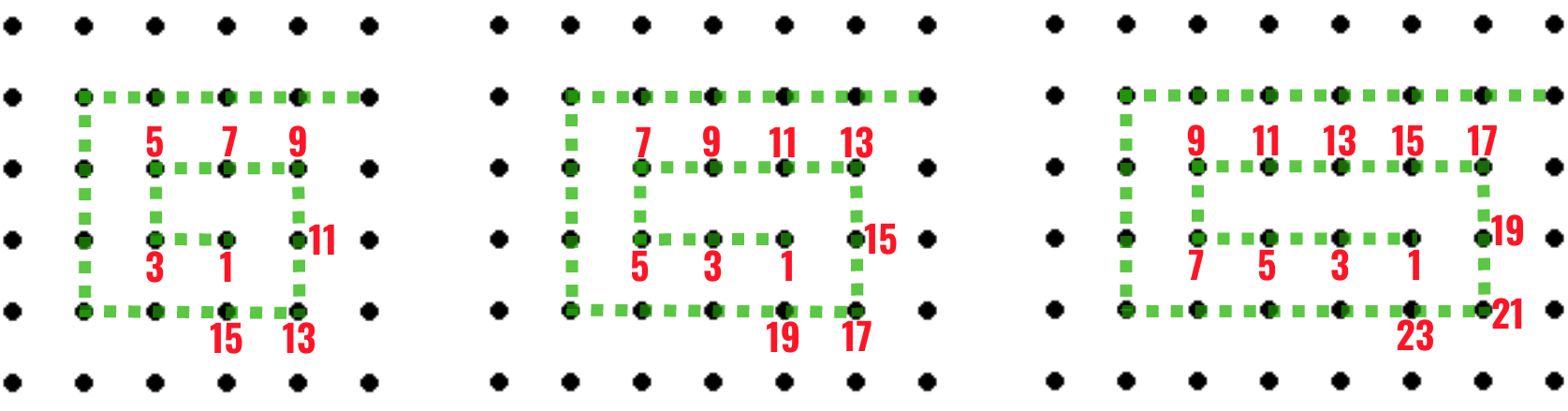}
\caption{The rectangular spiral and its development \protect \cite{2} for the cases of (from left to right) $n_3-n_2=0$, $n_3-n_2=1$, and $n_3-n_2=2$.
}
\label{fig:Figure_3S}
\end{figure}

It is possible to extend the aforementioned result we have previously shown to the general, $k$-dimensional, case $n_1 \times n_2 \times \cdots \times n_k$. The method to determine an acceptable upper limit for the optimal solution remains the same as in the case $n_1=n_2=\cdots =n_k$ so that
\begin{equation} \label{eq7}
h_u = (t+ 1) \cdot \prod_{j=1}^{k-3} n_j - 1,
\end{equation}
where $t$, the best upper bound available for the $n_{k-2} \times n_{k-1} \times n_k$ problem, is given by (\ref{eq4}) (except for the very particular cases we introduced at the beginning of Reference \cite{6}) and it is made explicit by (\ref{eq2},~\ref{eq6}).

Specifically, we will start considering an external grid defined by $[n_{k-1}; n_k]$, and we will connect the corresponding $n_k \cdot n_{k-1}$ points spending $2 \cdot n_{k-1}-1$ lines (following the rectangular spiral pattern); then, from the ending point of that external grid, we will draw the line segment which is orthogonal to any $[n_{k-1}; n_k]$ plane (along the $n_{k-2}$ points direction), and so forth.

The resulting rectangular spiral is a rectilinear spanning path that visits all the $n_1 \cdot n_2 \cdots n_k$ points of the grid $G(n_1, n_2, \ldots, n_k) \coloneqq \{ \{0,1,\ldots,n_1-1\} \times \{0,1,\ldots,n_2-1\} \times \cdots \times \{0,1,\ldots,n_k-1\}\} \subset \mathbb{R}^k$ only once and whose total Euclidean length is minimal among all the possible covering paths/trails of $G(n_1, n_2, \ldots, n_k)$ since it is given by $\prod_{i=1}^{k} n_i - 1$.

Consequently, we finally get a minimum length uncrossing covering path consisting of $(2 \cdot n_1 \cdot n_2 - c) \cdot \prod_{j=1}^{k-3} n_j -1$ axis-parallel line-segments (see (\ref{eq4})), entirely contained in the axis-aligned bounding box $[0,1]^k$, which provides a non-trivial upper bound for the general, unconstrained, $n_1 \times n_2 \times \cdots \times n_k$ Points Problem.


\section{The \texorpdfstring{$n_1 \times n_2 \times \cdots \times n_k$}{} Problem bounded from below}

In this section, we provide a lower bound for the ($k$-dimensional) $n_1 \times n_2 \times \cdots \times n_k$ Points Problem. In this way, we can build a range in which all the best possible solutions to the problem we are considering (for any given $k$-tuple of positive integers) will certainly fall.
In conclusion, we provide a few characteristic numerical examples in order to appreciate the quality of the result arising from the particular approach we have chosen.

For $k=3$ (and $n_1 \leq n_2 \leq n_3$, as usual), let us first examine the structure of the grid: it is not possible to intersect more than $(n_3-1)+(n_2-1)=n_3+n_2-2$ points using two consecutive lines; however, there is one exception (which, for the sake of simplicity, we may assume as in the case of the first two lines drawn). In this circumstance, it is possible to fit $n_3$ points with the first line and $n_2-1$ points using the second one, just as in the case of the pure rectangular spiral solution that we have already considered.

Now, let us observe that, lying (by definition) each segment on a unique plan, it will be necessary to provide $n_1-1$ lines to connect the various plans that are addressed in succession (of any type): there is no way to avoid using fewer than $n_1-1$ lines to connect (at most) $n_1-1$ points at a time (under the constraint previously explained above to connect $n_3+n_2-1$ points with the first two line segments). Each of these lines could be interposed between as many rectilinear line segments capable of connecting $n_k-1$  points at any one time.

Following the same pattern, we notice that the mentioned outcome does not substantially change in $4$, $5$, $6$, and more dimensions.

Let $h_l$ denote the number of line segments of our lower bound so that, for any $k \geq 3$, we have
\begin{equation} \label{eq8}
\prod_{i=1}^{k} n_i \coloneqq n_k+
\sum_{j=1}^{k-2} (n_j-1)^2 +(n_k-1) \cdot \sum_{j=1}^{k-2}(n_j-1) +\left[h_l-2 \cdot \sum_{j=1}^{k-2}(n_j-1)-1 \right] \cdot \floor[\bigg]{\frac{n_k+n_{k-1}}{2}-1}.
\end{equation}

Taking into account the fact that $\floor[\bigg]{\frac{n_k+n_{k-1}}{2}-1} \leq \floor[\bigg]{\frac{n_k+n_{k-1}-1}{2}}$, by doing some basic calculations, we get
\begin{equation}  \nonumber
\begin{aligned}
h_l =
    \begin{cases}
\ceil[\Bigg]{\frac{2 \cdot \left[\prod_{i=1}^{k} n_i -\sum_{j=1}^{k-2} {n_j}^2+(3-n_k ) \cdot \sum_{j=1}^{k-2}n_j+n_k \cdot (k-3)-2 \cdot k+4+(n_k+n_{k-1}-2) \cdot \sum_{j=1}^{k-2}(n_j-1) \right]}{n_k+n_{k-1}-2}}+1
\\ \mathbf{if} \hspace{5mm} \frac{n_k+n_{k-1}}{2} \in \mathbb{N}-\{0,1\}  \vspace{5mm}  \\
\ceil[\Bigg]{\frac{2 \cdot \left[\prod_{i=1}^{k} n_i -\sum_{j=1}^{k-2} {n_j}^2 +(3-n_k ) \cdot \sum_{j=1}^{k-2} n_j+n_k \cdot (k-3)-2 \cdot k+4+(n_k+n_{k-1}-1) \cdot \sum_{j=1}^{k-2}(n_j-1) \right]}{n_k+n_{k-1}-1}}+1 \\ \mathbf{if} \hspace{5mm} \frac{n_k+n_{k-1}+1}{2} \in \mathbb{N}-\{0,1\}
    \end{cases}.
\end{aligned}
\end{equation}

Hence,
\begin{equation}  \label{eq9}
\begin{aligned}
h_l =
    \begin{cases}
\ceil[\Bigg]{\frac{2 \cdot \left[\prod_{i=1}^{k} n_i -\sum_{j=1}^{k-2} {n_j}^2+\sum_{j=1}^{k-2}n_j- n_k +n_{k-1} \cdot \left(\sum_{j=1}^{k-2} n_j-k+2 \right) \right]}{n_k+n_{k-1}-2} }+1
\\ \mathbf{if} \hspace{5mm} \frac{n_k+n_{k-1}}{2} \in \mathbb{N}-\{0,1\}  \vspace{5mm}  \\
\ceil[\Bigg]{\frac{2 \cdot \left[\prod_{i=1}^{k} n_i -\sum_{j=1}^{k-2} {n_j}^2+2 \cdot \sum_{j=1}^{k-2} n_j- n_k +n_{k-1} \cdot \left(\sum_{j=1}^{k-2} n_j-k+2 \right)- k+2 \right]}{n_k+n_{k-1}-1} }+1 \\ \mathbf{if} \hspace{5mm} \frac{n_k+n_{k-1}+1}{2} \in \mathbb{N}-\{0,1\}
    \end{cases}.
\end{aligned}
\end{equation}

In detail (looking at (\ref{eq9})), if $k=3$ is given, it follows that
\begin{equation}  \label{eq10}
\begin{aligned}
h_l =
    \begin{cases}
\ceil[\Bigg]{\frac{2 \cdot  \left(n_1 \cdot  n_2 \cdot  n_3-{n_1}^2+n_1 \cdot n_2+n_1- n_2 -n_3 \right)}{n_3+n_2-2}}+1
\\ \mathbf{if} \hspace{5mm} \frac{n_3+n_2}{2} \in \mathbb{N}-\{0,1\}  \vspace{5mm}  \\
\ceil[\Bigg]{\frac{2 \cdot  \left(n_1 \cdot  n_2 \cdot  n_3-{n_1}^2+n_1 \cdot n_2+2 \cdot n_1-n_2-n_3 -1 \right)}{n_3+n_2-1}}+1 \\ \mathbf{if} \hspace{5mm} \frac{n_3+n_2+1}{2} \in \mathbb{N}-\{0,1\}
    \end{cases}.
\end{aligned}
\end{equation}

On specifics, for $n \coloneqq n_3 =n_2 = n_1 \geq 2$ (which implies $\frac{n_3+n_2}{2} \in \mathbb{N}-\{0,1\}$ since $\frac{n_3+n_2}{2} = n$), we get the trivial bound
\begin{equation}  \label{eq11}
h_l = n^2+n+1.   
\end{equation}

Then, as long as $k>2$ and $n_k \geq n_{k-1}\geq \cdots \geq n_2 \geq n_1 \geq 2$,
\begin{gather*}
\frac{2}{n_k+n_{k-1}-2} \cdot \left[\prod_{i=1}^{k} n_i -\sum_{j=1}^{k-2} {n_j}^2+\sum_{j=1}^{k-2} n_j-n_k +n_{k-1} \cdot \left(\sum_{j=1}^{k-2} n_j-k+2 \right) \right] \geq  \\
\frac{2}{n_k+n_{k-1}-1} \cdot \left[\prod_{i=1}^{k} n_i -\sum_{j=1}^{k-2} {n_j}^2+2 \cdot \sum_{j=1}^{k-2} n_j-n_k +n_{k-1} \cdot \left(\sum_{j=1}^{k-2} n_j-k+2 \right)-k+2 \right].
\end{gather*}

Consequently, without loss of generality,
\begin{equation}  \label{eq12}
h_l \geq \ceil[\Bigg]{\frac{2}{n_k+n_{k-1}-1} \cdot \Biggl[\prod_{i=1}^{k} n_i -\sum_{j=1}^{k-2}{n_j}^2+2 \cdot \sum_{j=1}^{k-2} n_j-n_k +n_{k-1} \cdot \Biggl(\sum_{j=1}^{k-2} n_j-k+2 \Biggr)-k+2 \Biggr] }+1
\end{equation}
holds for each $k$-tuple $(n_k,n_{k-1},\ldots,n_2,n_1)$.


\section{Conclusion}

Given $k=3$, by combining (\ref{eq12}) and (\ref{eq2}, \ref{eq6}), we get the intervals in which the best possible solutions to the problem will fall (and let us point out that the lower bound from (\ref{eq12}) can also be improved by using (\ref{eq9})).

The width of the resulting range (and consequently how good our outcome may be) depends on the selected $3$-tuple $(n_1,n_2,n_3)$. Additionally, for any given pair $(n_1,n_2)$, our stated lower bound will always match the rectangular spiral upper bound as long as $n_3 \coloneqq n_3(n_1,n_2)$ is sufficiently large.

Example 1: $n_1=10$, $n_2=13$, $n_3=15$. Then,
$$147 \leq h \leq 253.$$

Example 2: $n_1=10$, $n_2=21$, $n_3=174$. Then,
$$378 \leq h \leq 419.$$

If $k>3$, the interval is given by
\begin{gather*}
\ceil[\Bigg]{\frac{2 \cdot \left(\prod_{i=1}^{k} n_i -\sum_{j=1}^{k-2}  {n_j}^2+2 \cdot \sum_{j=1}^{k-2} n_j- n_k +n_{k-1} \cdot \left(\sum_{j=1}^{k-2}  n_j-k+2 \right)-k+2 \right)}{n_k+n_{k-1}-1} }+1 \leq \\ h \leq (t+1) \cdot \prod_{j=1}^{k-3} n_j -1,
\end{gather*}
where $t$, the minimum upper limit for the $n_{k-2} \times n_{k-1} \times n_k$ Points Problem, is the result obtained from (\ref{eq4},~\ref{eq6}).

In this case, the size of the resulting interval also depends on the particular value of $k$ (generally speaking, the larger the $k$, the wider the interval).

Example 3: $k=4$, $n_1=10$, $n_2=16$, $n_3=18$, $n_4=48$ (and thus $t=575$). Then,
$$4257 \leq h \leq 5759.$$

Furthermore, for some particular $k$-tuples $(n_1,n_2,\ldots,n_k)$, the upper and the lower bound coincide, thus allowing us to obtain a complete and definitive resolution of the given problem.\newline
E.g., if $k=3$, $n_1=n_2=3$, and $n_3 \geq 27$, then $h_l=h_u=h_{best}=17$. Ditto if $k=3$, $n_1=3$, $n_2=4$, and $n_3 \geq 56$, as $h_l=h_u=h_{best}=23$. While, if $k=4$, $n_1=n_2=n_3=2$, and $n_4 \geq 10$, $h_l=h_u=h_{best}=15$ follows.
\vspace{3mm}

\noindent [We omit the last pages of the 2014 version of this paper, the appendix and the improved bounds for special cases only since better results are now available (e.g., as long as $k \geq 3$ and $n_1=n_2=\cdots=n_k$, Equation (6) of Reference \cite{11} provides the upper bound $\left(\floor*{\frac{3}{2} \cdot {n_1}^2} - \floor*{\frac{n_1-1}{4}} + \floor*{\frac{n_1+1}{4}} - \floor*{\frac{n_1+2}{4}} + \floor*{\frac{n_1}{4}} + n_1 - 1 \right) \cdot {n_1}^{k-3} - 1$ while Equations~(4) and (5) state also that $3 \cdot 2^{k-2}$ \cite{12} and $\frac{3^k-1}{2}$ \cite{8} are the exact solutions for $2 = n_1=n_2=\cdots=n_k$ and $3 = n_1=n_2=\cdots=n_k$, respectively).]


\section*{Acknowledgments}

The author sincerely thanks Graham Powell, MA, for his constant support.

\makeatletter
\renewcommand{\@biblabel}[1]{[#1]\hfill}

\end{document}